\documentclass[11pt,reqno]{amsart}
\usepackage{color}

\usepackage{amsmath}
\usepackage{amsfonts,amscd}
\usepackage{amssymb}
\usepackage{url}
\usepackage{hyperref}

\usepackage[english]{babel}
\usepackage{booktabs}
\usepackage{tikz}

\theoremstyle{plain}
\newtheorem{theorem}                 {Theorem}      [section]

\theoremstyle{definition}

\newtheorem{definition}   [theorem]  {Definition}

\newtheorem{example}      [theorem]  {Example}

\numberwithin{equation}{section}

\def \cn{{\mathbb C}}

\def \hn{{\mathbb H}}

\def \rn{{\mathbb R}}

\def \zn{{\mathbb Z}}

\def \E{\mathcal E}

\def \P{\mathcal P}

\def\nab#1#2{\hbox{$\nabla$\kern -.3em\lower 1.0 ex
		\hbox{$#1$}\kern -.1 em {$#2$}}}
\def\hatnab#1#2{\hbox{$\nabla$\kern -.3em\lower 1.0 ex
		\hbox{$#1$}\kern -.1 em {$#2$}}}

\def \Re{\mathfrak R\mathfrak e}

\def \ip #1#2{\langle #1,#2 \rangle}
\def \ipp #1#2{\langle\hskip -.1cm\langle #1,#2 \rangle\hskip -.1cm\rangle}

\def \SL2{\widetilde{\text{\bf SL}}_{2}(\rn)}

\def \SO#1{\mathbf{SO}(#1)}

\def \U#1{\text{\bf U}(#1)}

\def \SU#1{\text{\bf SU}(#1)}

\def \Sp#1{\text{\bf Sp}(#1)}

\DeclareMathOperator{\Div}{div}

\numberwithin{equation}{section}
\allowdisplaybreaks

\begin{document}

	

\author{Sigmundur Gudmundsson}
\address{Mathematics, Faculty of Science\\
Lund University\\
Box 118, Lund 221 00\\
Sweden}
\email{Sigmundur.Gudmundsson@math.lu.se}

\title
[Minimal Submanifolds]
{Minimal Submanifolds of The Complex and Quaternionic Projective and Hyperbolic Spaces $\cn P^{2n-1}$, $\hn P^{2n-1}$, $\cn H^{2n-1}$, $\hn H^{2n-1}$\\ via Harmonic Morphisms}

\begin{abstract}
In this work we construct non-holomorphic, complete and minimal submanifolds of the odd-dimensional complex projective spaces $\cn P^{2n-1}$ and their dual complex hyperbolic spaces $\cn H^{2n-1}$.  We then provide complete minimal submanifolds of the quaternionic projective spaces $\hn P^{2n-1}$ and their dual quaternionic hyperbolic spaces $\hn H^{2n-1}$. All the constructed minimal submanifolds are of codimension two. Our main tools are complex-valued harmonic morphisms from the above mentioned ambient spaces.
\end{abstract}
	
{\Large }\dedicatory{version 1.004 - \today}

\maketitle

\section{Introduction}
\label{section-introduction}

The study of minimal submanifolds of a given ambient space plays a central role in differential geometry.  This has a long, interesting history and has attracted the interests of profound mathematicians for many generations.  The famous Weierstrass-Enneper representation formula, for minimal surfaces in three-dimensional Euclidean space, brings {\it complex analysis} into play as a useful tool for the study of these beautiful objects.

This was later generalised to the study of minimal surfaces in much more general ambient manifolds via {\it harmonic conformal immersions}.  The next  result follows from the seminal paper \cite{Eel-Sam} of Eells and Sampson from 1964.  For this see also Proposition 3.5.1 of \cite{Bai-Woo-book}.

\begin{theorem}
	Let $\phi:(M^m,g)\to (N,h)$ be a smooth conformal map between Riemannian manifolds.  If $m=2$ then $\phi$ is harmonic if and only if the image is minimal in $(N,h)$.
\end{theorem}

This result has turned out to be very useful in the construction of minimal surfaces in Riemannian symmetric spaces of various types.  For this we refer to \cite{Cal},
\cite{Eel-Woo}, \cite{Uhl}, \cite{Bur-Raw} and \cite{Bur-Gue}, just to name a few.
\smallskip

In their work \cite{Bai-Eel} from 1981, Baird and Eells have shown that complex-valued harmonic morphisms from Riemannian manifolds are useful tools for the study of minimal submanifolds of codimension two.  The following result is a slightly adapted version to the semi-Riemannian situation.

\begin{theorem}\label{theorem-Bai-Eel-special}\cite{Gud-8}
Let $\phi:(M,g)\to\cn$ be a complex-valued harmonic morphism from a semi-Riemannian manifold.  Then every regular fibre of $\phi$ is a minimal submanifold of $(M,g)$ of codimension {\it two}.
\end{theorem}

This can be seen as dual to the above-mentioned generalisation of the Weierstrass-Enneper representation. Harmonic morphisms are the much studied {\it horizontally conformal harmonic maps}.  For an introduction to the general theory we recommened the book \cite{Bai-Woo-book}, by Baird and Wood, and the regularly updated online bibliography \cite{Gud-bib}.
\medskip

In this work we construct minimal submanifolds of the odd-dimensional complex projective spaces $\cn P^{2n-1}$, with  $n\ge 2$.  These are compact and non-holomorphic submanifolds of the K\" ahler manifold $\cn  P^{2n-1}$, see Example \ref{example-complex-projective}.  The complex projective spaces $\cn P^{2n-1}$ are well-known Riemannian symmetric spaces of rank one and so are their complex hyperbolic dual spaces $\cn H^{2n-1}$, see \cite{Hel}.  For such a situation we have a duality principle for complex-valued harmonic morphisms, developed in \cite{Gud-Sve-1}, see also \cite {Gud-Mon-Rat-1}.  This means that we can translate the solutions found in the compact projective cases directly over to the non-compact hyperbolic spaces.  With this we yield complete, minimal and non-holomorphic submanifolds of the complex K\" ahler manifolds $\cn H^{2n-1}$, with $n\ge 2$.  For this see Example \ref{example-complex-hyperbolic}.

We then turn our focus to the corresponding quaternionic Riemannian symmetric spaces of rank one.  We construct  compact minimal submanifolds of the quaternionic projective spaces $\hn P^{n-1}$, see Example \ref{example-quaternionic-projective}.  We then apply the duality principle to provide complete minimal submanifolds of the quaternionic hyperbolic spaces $\hn H^{n-1}$, see Example \ref{example-quaternionic-hyperbolic}

All the constructed minimal submanifolds are of codimension two.  Our most important tools are the complex-valued harmonic morphisms and the above mentioned Theorem \ref{theorem-Bai-Eel-special}.

\section{Preliminaries}

Let $M$ and $N$ be two manifolds of dimensions $m$ and $n$,
respectively. Then a semi-Riemannian metric $g$ on $M$ gives rise
to the notion of a Laplace-Beltrami operator (alt. tension field) on $(M,g)$ and real-valued harmonic
functions $f:(M,g)\to\rn$. This can be generalised to the concept
of a harmonic map $\phi:(M,g)\to (N,h)$ between semi-Riemannian
manifolds being a solution to a non-linear system of partial
differential equations. 

\begin{definition}
A map $\phi:(M,g)\to (N,h)$ between semi-Riemannian manifolds is called a {\it harmonic morphism} if, for any harmonic function $f:U\to\rn$ defined on an open subset $U$ of $N$ with $\phi^{-1}(U)$ non-empty, the composition $f\circ\phi:\phi^{-1}(U)\to\rn$ is a harmonic function.
\end{definition}

The following characterisation of the harmonic morphisms between semi- Riemannian manifolds is due to Fuglede, see \cite{Fug-2}. For the definition of horizontal conformality we
refer to \cite{Bai-Woo-book}.

\begin{theorem}\cite{Fug-2}
A map $\phi:(M,g)\to (N,h)$ between semi-Riemannian manifolds is a harmonic morphism if and only if it is a horizontally (weakly) conformal harmonic map.
\end{theorem}

The next result generalises the corresponding well-known theorem
of Baird and Eells in the Riemannian case, see \cite{Bai-Eel}. It
gives the theory of harmonic morphisms a strong geometric flavour
and shows that the case when $n=2$ is particularly interesting. In
that case the conditions characterising harmonic morphisms are
independent of conformal changes of the metric on the surface
$N^2$.  For the definition of horizontal homothety we refer to
\cite{Bai-Woo-book}.

\begin{theorem}\cite{Gud-8}\label{theorem-minimal}
Let $(M^m,g)$ be a semi-Riemannian manifold, $(N^n,h)$ be Riemannian and $\pi:M\to N$ be a horizontally conformal submersion.  If
\begin{enumerate}
\item[i.] $n=2$, then $\pi$ is harmonic if and only if $\pi$ has minimal fibres,
\item[ii.] $n\ge 3$, then two of the following conditions imply the other,
\begin{enumerate}
\item $\pi$ is a harmonic map,
\item $\pi$ has minimal fibres, 
\item $\pi$ is horizontally homothetic.	
\end{enumerate}
\end{enumerate}
\end{theorem}

The next result shows that an eigenfamily on a semi-Riemannian manifold,  see Section \ref{section-eigenfunctions}, can be used to produce a variety of local harmonic morphisms.

\begin{theorem}\cite{Gud-Sak-1}\label{theorem-rational}
Let $(M,g)$ be a semi-Riemannian manifold and 
$$\E =\{\phi_1,\dots,\phi_n\}$$ 
be a finite eigenfamily of complex-valued functions on $M$. If $P,Q:\cn^n\to\cn$ are linearily independent homogeneous polynomials of the same positive degree then the quotient
$$\frac{P(\phi_1,\dots ,\phi_n)}{Q(\phi_1,\dots ,\phi_n)}$$ 
is a non-constant harmonic morphism on the open and dense subset	$$\Omega (Q)=\{p\in M\,|\, Q(\phi_1(p),\dots ,\phi_n(p))\neq 0\}.$$
\end{theorem}

\section{Eigenfunctions and Eigenfamilies}
\label{section-eigenfunctions}

Let $(M,g)$ be an $m$-dimensional semi-Riemannian manifold and $T^{\cn}M$ be the complexification of the tangent bundle $TM$ of $M$. We extend the metric $g$ to a complex bilinear form on $T^{\cn}M$.  Then the gradient $\nabla\phi$ of a complex-valued function $\phi:(M,g)\to\cn$ is a section of $T^{\cn}M$.  In this situation, we have the well-known complex linear {\it Laplace-Beltrami operator} (alt. {\it tension field}) $\tau$ on $(M,g)$.  In local coordinates this satisfies 
$$
\tau(\phi)=\Div (\nabla \phi)=\sum_{i,j=1}^m\frac{1}{\sqrt{|g|}} \frac{\partial}{\partial x_j}
\left(g^{ij}\, \sqrt{|g|}\, \frac{\partial \phi}{\partial x_i}\right).
$$
For two complex-valued functions $\phi,\psi:(M,g)\to\cn$ we have the following well-known fundamental relation
\begin{equation*}\label{equation-basic}
	\tau(\phi\cdot \psi)=\tau(\phi)\cdot\psi +2\,\kappa(\phi,\psi)+\phi\cdot\tau(\psi),
\end{equation*}
where the symmetric complex bilinear {\it conformality operator} $\kappa$ is given by $$\kappa(\phi,\psi)=g(\nabla \phi,\nabla \psi).$$  Locally this satisfies 
$$\kappa(\phi,\psi)=\sum_{i,j=1}^mg^{ij}\cdot\frac{\partial\phi}{\partial x_i}\frac{\partial \psi}{\partial x_j}.$$

The naming of the operator $\kappa$ comes from the fact that $\kappa (\phi,\phi)=0$ if and only if 
$$\kappa (\phi,\phi)=|\nabla u|^2-|\nabla v|^2+2i\cdot  g(\nabla u,\nabla v)=0.$$

\begin{definition}\cite{Gud-Sak-1}\label{definition-eigenfamily}
Let $(M,g)$ be a semi-Riemannian manifold. Then a complex-valued function $\phi:M\to\cn$ is said to be a {\it $(\lambda,\mu)$-eigenfunction} if it is eigen both with respect to the Laplace-Beltrami operator $\tau$ and the conformality operator $\kappa$ i.e. there exist complex numbers $\lambda,\mu\in\cn$ such that $$\tau(\phi)=\lambda\cdot\phi\ \ \text{and}\ \ \kappa(\phi,\phi)=\mu\cdot \phi^2.$$	
A set $\E =\{\phi_i:M\to\cn\ |\ i\in I\}$ of complex-valued functions is said to be a {\it $(\lambda,\mu)$-eigenfamily} on $M$ if there exist complex numbers $\lambda,\mu\in\cn$ such that for all $\phi,\psi\in\E$ we have 
$$\tau(\phi)=\lambda\cdot\phi\ \ \text{and}\ \ \kappa(\phi,\psi)=\mu\cdot \phi\,\psi.$$ 
\end{definition}

For the standard odd-dimensional round spheres we have the following eigenfamilies based on the classical real-valued spherical harmonics.

\begin{example}\cite{Gud-Mun-1}
\label{example-basic-sphere} 
Let $S^{2n-1}$ be the odd-dimensional unit sphere in the standard Euclidean space $\cn^{n}\cong\rn^{2n}$ and define $\phi_1,\dots,\phi_n:S^{2n-1}\to\cn$ by
$$\phi_j:(z_1,\dots,z_{n})\mapsto \frac{z_j}{\sqrt{|z_1|^2+\cdots +|z_n|^2}}.$$  Then the tension field $\tau$ and the conformality operator $\kappa$ on $S^{2n-1}$ satisfy	$$\tau(\phi_j)=-\,(2n-1)\cdot\phi_j\ \ \text{and}\ \ \kappa(\phi_j,\phi_k)=-\,1\cdot \phi_j\cdot\phi_k.$$
\end{example}

With the following result we show that a given eigenfamily $\E$ can be used to produce a large collection $\P_d(\E)$ of such objects.

\begin{theorem}\label{theorem-polynomials}\cite{Gha-Gud-5}
Let $(M,g)$ be a semi-Riemannian manifold and the set of complex-valued functions  $$\E=\{\phi_i:M\to\cn\,|\,i=1,2,\dots,n\}$$ 
be a finite eigenfamily i.e. there exist complex numbers $\lambda,\mu\in\cn$ such that for all $\phi,\psi\in\E$ $$\tau(\phi)=\lambda\cdot\phi\ \ \text{and}\ \ \kappa(\phi,\psi)=\mu\cdot\phi\,\psi.$$  
Then the set of complex homogeneous polynomials of degree $d$
$$\P_d(\E)=\{P:M\to\cn\,|\, P\in\cn[\phi_1,\phi_2,\dots,\phi_n],\, P(\alpha\cdot\phi)=\alpha^d\cdot P(\phi),\, \alpha\in\cn\}$$ 
is an eigenfamily on $M$ such that for all $P,Q\in\P_d(\E)$ we have
$$\tau(P)=(d\,\lambda+d(d-1)\,\mu)\cdot P\ \ \text{and}\ \ \kappa(P,Q)=d^2\mu\cdot P\, Q.$$
\end{theorem}

\section{The Euclidean Space $\cn^{2n}$}

Let $\cn^{2n}$ be the complex $2n$-dimensional vector space equipped with its standard Euclidean metric $\ip \cdot\cdot :\cn^{2n}\times\cn^{2n}\to\rn$ satisfying 
$$\ip zw=\sum_{k=1}^{2n}\Re (z_k\cdot\bar w_k).$$
For two complex-valued functions $\hat\phi,\hat\psi:\cn^{2n}\to\cn$ the tension field $\tau$ and the conformal operator $\kappa$ are given by 
$$\tau(\hat\phi)=4\cdot\sum_{k=1}^{2n}\frac{\partial^2 \hat\phi}{\partial z_k\partial\bar z_k}
\ \ \text{and}\ \ \kappa(\hat\phi,\hat\psi)=2\cdot\sum_{k=1}^{2n}\big(\frac{\partial\hat\phi}{\partial z_k}\frac{\partial\hat\psi}{\partial\bar z_k}+\frac{\partial\hat\phi}{\partial \bar z_k}\frac{\partial\hat\psi}{\partial z_k}\big).$$
Here $z=(z_1,z_2,\dots ,z_{2n})$ are the standard global coordinates on $\cn^{2n}$.

\begin{example}\cite{Gud-Mun-1}
\label{example-basic-complex-space}
For two integers $j,k\in\zn$ satisfying $1\le j\le n$ and  $n+1\le k\le 2n$, we define the function $\hat\phi_{jk}:\cn^{2n}\to\cn$ with 
$$
\hat\phi_{jk}:(z_1,\dots,z_{2n})\mapsto z_j\cdot\bar z_k.$$ 
Then it is easily seen that we have a complex $n^2$-dimensional $(0,0)$-eigenfamily $\hat\E_n$ on $\cn^{2n}$ satisfying 
$$
\hat\E_n=\{\hat\phi_{jk}\,|\,1\le j\le n\ \ \text{and}\ \ n+1\le k\le 2n\}.
$$
\end{example}

We will now illustrate our new construction method in the special case when $n=3$. It is easily seen that it actually  works for any $n\ge 2$.

\begin{example}\label{example-complex-projective-main}
Let $d\in\zn^+$ be an arbitrary positive integer and $a,b\in\cn^9$ be two linearly independent elements, satisfying
$$
a=(a_{14},a_{15},a_{16},a_{24},a_{25},a_{26},a_{34},a_{35},a_{36}),
$$
$$
b=(
b_{14},b_{15},b_{16},b_{24},b_{25},b_{26},b_{34},b_{35},b_{36}).
$$
Further let the matrices $A,B\in\cn^{3\times 3}$ be given by 
$$A=
\begin{bmatrix}
a_{14} & a_{15} & a_{16}  \\
a_{24} & a_{25} & a_{26}  \\
a_{34} & a_{35} & a_{36}
\end{bmatrix},\ \ B=
\begin{bmatrix}
b_{14} & b_{15} & b_{16}  \\
b_{24} & b_{25} & b_{26}  \\
b_{34} & b_{35} & b_{36}
\end{bmatrix}
$$
and assume that $\det B\neq 0$.  Then the third order polynomial $R:\cn\to\cn$ with $R(s)=\det (s\cdot B-A)$ is of the form 
$$R(s)=\det B\cdot s^3-c_2\cdot s^2+c_1\cdot s-\det A,$$ for some complex numbers $c_1,c_2\in\cn$.
\smallskip

For $j=1,2,3$ and $k=4,5,6$ let us now define the two linearly independent polynomials $\hat P,\hat Q:\cn^9\to\cn$ by
$$
\hat P(z)=\sum_{j,k}a_{jk}\cdot (z_j\bar z_k)^d\ \ \text{and}\ \ 
\hat Q(z)=\sum_{j,k}b_{jk}\cdot (z_j\bar z_k)^d.
$$
Further we define the complex-valued function $\hat F:\Omega (\hat Q)\to\cn$ by $\hat F(z)=\hat P(z)/\hat Q(z)$, where $\Omega(\hat Q)=\{z\in\cn^6\, |\, \hat Q(z)\neq 0\}$.
\smallskip

Let $\alpha\in\cn^*$ be a non-zero complex number such that $R(\alpha)\neq 0$.  Then a point $z\in \hat F^{-1}(\{\alpha\})$ in the fibre over $\alpha$ is critical if and only if 
$$\frac{\partial\hat F}{\partial z_1}(z)
=\frac{\partial\hat F}{\partial z_2}(z)
=\frac{\partial\hat F}{\partial z_3}(z)=0
=\frac{\partial\hat F}{\partial \bar z_4}(z)
=\frac{\partial\hat F}{\partial \bar z_5}(z)
=\frac{\partial\hat F}{\partial \bar z_6}(z).$$
It is easily seen that this is equivalent to the following systems of equations
$$
\Big( t\cdot B-A\Big)\cdot
\begin{bmatrix}
	z_1^{d-1}\bar z_4^d & z_1^d\bar z_4^{d-1} \\[.5em]
	z_2^{d-1}\bar z_5^d & z_2^d\bar z_5^{d-1} \\[.5em]
	z_3^{d-1}\bar z_6^d & z_3^d\bar z_6^{d-1}
\end{bmatrix}=0.
$$

Since $R(\alpha)\neq 0$ this implies that $\hat Q(z)=0$ so $z$ is not contained in $\Omega (\hat Q)$. This shows that every point $z\in \hat F^{-1}(\{\alpha\})$ is regular.  It now follows from Theorem \ref{theorem-minimal} and Theorem \ref{theorem-rational} that the fibre 
$$
\hat F^{-1}(\{\alpha\})=\{z\in\Omega(\hat F)\,|\, \hat P(z)=\alpha\cdot\hat Q(z)\},
$$ 
over the non-zero element $\alpha\in\cn^*$, is a complete {\it minimal submanifold} of $\cn^6$ of codimension two.
\end{example}
\smallskip


\section{The Complex Projective Space $\cn P^{2n-1}$}

In this section we construct compact non-holomorphic minimal submanifolds, of the odd-dimensional complex projective space $\cn P^{2n-1}$, of codimension two. Here we use the construction presented in Example \ref{example-complex-projective-main}.
\medskip 

The map $\rho:S^1\times S^{4n-1}\to S^{4n-1}$ with $\rho:(e^{i\theta},z)\mapsto e^{i\theta} z$ is a smooth action of the Lie group $S^1$ on the unit sphere $S^{4n-1}$ in $\cn^{2n}$.  The quotient space of this action is the well-known complex projective space $\cn P^{2n-1}$.  The natural projection $\pi:S^{4n-1}\to\cn P^{2n-1}$ is a Riemannian submersion with totally geodesic fibres and hence a harmonic morphism.

\begin{example}\label{example-complex-projective}
For $j,k\in\zn$ satisfying $1\le j\le n$ and  $n+1\le k\le 2n$, we define the $S^1$-invariant function $\tilde\phi_{jk}:S^{4n-1}\to\cn$ by 
$$
\tilde\phi_{jk}:(z_1,\dots,z_{2n})\mapsto \frac
{z_j\cdot\bar z_k}{\ip zz}.
$$  
Then $\tilde\phi_{jk}:S^{4n-1}\to\cn$ is the restriction of 
the function $\hat\phi_{jk}:\cn^{2n}\to\cn$ to the unit sphere $S^{4n-1}$. An easy calculation shows that the tension field $\tau$ and the conformality operator $\kappa$ on $S^{4n-1}$ satisfy
$$
\tau(\tilde\phi_{jk})=-\,8n\cdot\tilde\phi_{jk}\ \ \text{and}\ \ \kappa(\tilde\phi_{jk},\tilde\phi_{lm})=-\,4\cdot \tilde\phi_{jk}\,\tilde\phi_{lm}.
$$
Hence we have a complex $n^2$-dimensional eigenfamily $\tilde\E_n$ on the unit sphere $S^{4n-1}$. 

This implies that the restrictions $\tilde P,\tilde Q:S^{4n-1}\to\cn$ of the polynomials $\hat P,\hat Q:\cn^{2n}\to\cn$ to the unit sphere $S^{4n-1}$ form an eigenfamily.  Then we define the complex-valued harmonic morphism $\tilde F:\Omega (\tilde Q)\to\cn$ by $\tilde F(z)=\tilde P(z)/\tilde Q(z)$, where 
$$
\Omega(\tilde Q)=\{z\in S^{4n-1}\, |\, \tilde Q(z)\neq 0\}.
$$

The polynomial maps $\hat P,\hat Q:\cn^{2n}\to\cn$ on $\cn^{2n}$ induce  linearly independent functions $P,Q:\cn P^{2n-1}\to\cn$ on the complex projective space $\cn P^{2n-1}$, satisfying
$$
P:[z]\mapsto\sum_{j,k}a_{jk}\cdot\phi_{jk}^d\ \ \text{and}\ \ 
Q:[z]\mapsto\sum_{j,k}b_{jk}\cdot\phi_{jk}^d.
$$
Here the coefficients $a_{jk},b_{jk}\in\cn$ are chosen as explained in Example \ref{example-complex-projective} and the induced functions $\phi_{jk}:\cn P^{2n-1}\to\cn$ satisfy 
$$
\phi_{jk}:[z]\mapsto \frac
{z_j\cdot\bar z_k}{\ip zz}.
$$ 
The natural projection $\pi:S^{4n-1}\to\cn P^{2n-1}$ is a harmonic morphism so the complex valued functions $P,Q$ form an eigenfamily $\E_n$ on the complex projective space $\cn P^{2n-1}$.
\smallskip 

Let $\Omega (Q)$ be the subset of $\cn P^{2n-1}$ given by 
$$
\Omega (Q)=\{[z]\in\cn P^{2n-1}\,|\,Q([z])\neq 0\}
$$
and define $F:\Omega(Q)\to\cn$ by $F([z])=P([z])/Q([z])$.  Then for a non-zero element  $\alpha\in\cn^*$ with $R(\alpha)\neq 0$, the fibre $F^{-1}(\{\alpha\})$ satisfies $F^{-1}(\{\alpha\})=\pi (\tilde F^{-1}(\{\alpha\}))$.  Since $\pi:S^{4n-1}\to\cn P^{2n-1}$ is a Riemannian submersion the fibre is a {\it minimal submanifold} of codimension two satisfying 
$$
F^{-1}(\{\alpha\})=\{[z]\in\cn P^{2n-1}\,|\, P([z])=\alpha\cdot Q([z])\}.
$$
This is clearly a {\it compact} and {\it non-holomorpic} submanifold of $\cn P^{2n-1}$.
\end{example}

\section{The Semi-Euclidean Space $\cn^{2n}_1$}

Let $\cn^{2n}_1$ be the complex $2n$-dimensional vector space equipped with its standard semi-Euclidean metric $\ipp \cdot\cdot :\cn^{2n}_1\times\cn^{2n}_1\to\rn$ satisfying 
$$
\ipp zw=\Re (-z_1\cdot\bar w_1+\sum_{k=2}^{2n}z_k\cdot\bar w_k).
$$
Let $\Omega^{2n}_1$ be the open subset $\Omega^{2n}_1=\{z\in\cn^{2n}_1\,|\, \ipp zz<0\}$ of $\cn^{2n}_1$.
For two complex-valued functions $\hat\phi^*,\hat\psi^*:\Omega^{2n}_1\to\cn$ the tension field $\tau$ and the conformal operator $\kappa$ on $\Omega^{2n}_1$, are given by 
$$
\tau(\hat\phi^*)=4\cdot\big(-\frac{\partial^2 \hat\phi*}{\partial z_1\partial\bar z_1}+ \sum_{k=2}^{2n}\frac{\partial^2 \hat\phi^*}{\partial z_k\partial\bar z_k}\big),
$$
$$
\kappa(\hat\phi^*,\hat\psi^*)=2\cdot\Big(-\frac{\partial\hat\phi^*}{\partial z_1}\frac{\partial\hat\psi^*}{\partial\bar z_1}-\frac{\partial\hat\phi^*}{\partial \bar z_1}\frac{\partial\hat\psi^*}{\partial z_1}+\sum_{k=2}^{2n}\big(\frac{\partial\hat\phi^*}{\partial z_k}\frac{\partial\hat\psi^*}{\partial\bar z_k}+\frac{\partial\hat\phi^*}{\partial \bar z_k}\frac{\partial\hat\psi^*}{\partial z_k}\big)\Big).$$
Here $z=(z_1,z_2,\dots ,z_{2n})$ are the standard global coordinates on $\Omega^{2n}_1$.

\begin{example}\cite{Gud-Mun-1}\label{example-basic-complex-space}
For $j,k\in\zn$ satisfying $1\le j\le n$ and  $n+1\le k\le 2n$, define the functions  $\hat\phi^*_{jk}:\Omega^{2n}_1\to\cn$ with 
$$
\hat\phi^*_{jk}:(z_1,\dots,z_{n+1})\mapsto {z_j\cdot\bar z_k}.
$$ 
Then it is easily seen that we have a complex $n^2$-dimensional $(0,0)$-eigenfamily $\hat\E_n^*$ on $\cn^{2n}_1$ satisfying 
$$
\hat\E_n^*=\{\hat\phi_{jk}^*\,|\,1\le j\le n\ \ \text{and}\ \ n+1\le k\le 2n\}.
$$
\end{example}

\begin{example}\label{example-hyperbolic-main}
Let $d\in\zn^+$ be an arbitrary positive integer and $a,b$ be two linearly independent elements of the standard complex vector space $\cn^9$, satisfying
$$
a=(a_{14},a_{15},a_{16},a_{24},a_{25},a_{26},a_{34},a_{35},a_{36}),$$ 
$$
b=(b_{14},b_{15},b_{16},b_{24},b_{25},b_{26},b_{34},b_{35},b_{36}).
$$
Further let the matrices $A,B\in\cn^{3\times 3}$ be given by 
$$A=
\begin{bmatrix}
a_{14} & a_{15} & a_{16}  \\
a_{24} & a_{25} & a_{26}  \\
a_{34} & a_{35} & a_{36}
\end{bmatrix},\ \ B=
\begin{bmatrix}
b_{14} & b_{15} & b_{16}  \\
b_{24} & b_{25} & b_{26}  \\
b_{34} & b_{35} & b_{36}
\end{bmatrix}
$$
and assume that $\det B\neq 0$.  Then the third order polynomial $R:\cn\to\cn$ with $R(s)=\det (s\cdot B-A)$ is of the form 
$$
R(s)=\det B\cdot s^3-c_2\cdot s^2+c_1\cdot s-\det A,
$$ 
for some complex numbers $c_1,c_2\in\cn$.
\smallskip
	
For $j=1,2,3$ and $k=4,5,6$ let us now define the two linearly independent polynomials $\hat P^*,\hat Q^*:\Omega^6_1\to\cn$ by
$$
\hat P^*(z)=\sum_{j,k}a_{jk}\cdot (z_j\bar z_k)^d\ \ \text{and}\ \ 
\hat Q^*(z)=\sum_{j,k}b_{jk}\cdot (z_j\bar z_k)^d.
$$
Further we define the complex-valued function $\hat F^*:\Omega (\hat Q)\to\cn$ by $\hat F^*(z)=\hat P^*(z)/\hat Q^*(z)$, where $\Omega(\hat Q^*)=\{z\in\Omega^6_1\, |\, \hat Q^*(z)\neq 0\}$.
\smallskip
	
Let $\alpha\in\cn^*$ be a non-zero complex number such that $R(\alpha)\neq 0$.  Then a point $z\in (\hat F^*)^{-1}(\{\alpha\})$ in the fibre over $\alpha$ is critical if and only if 
$$
 \frac{\partial F^*}{\partial z_1}(z)
=\frac{\partial F^*}{\partial z_2}(z)
=\frac{\partial F^*}{\partial z_3}(z)=0
=\frac{\partial F^*}{\partial \bar z_4}(z)
=\frac{\partial F^*}{\partial \bar z_5}(z)
=\frac{\partial F^*}{\partial \bar z_6}(z).
$$
It is easily seen that this is equivalent to the following systems of equations

$$
\Big(t\cdot B-A\Big)\cdot
\begin{bmatrix}
z_1^{d-1}\bar z_4^d & z_1^d\bar z_4^{d-1} \\[.5em]
z_2^{d-1}\bar z_5^d & z_2^d\bar z_5^{d-1} \\[.5em]
z_3^{d-1}\bar z_6^d & z_3^d\bar z_6^{d-1}
\end{bmatrix}=0.
$$

Since $R(\alpha)\neq 0$ this implies that $\hat Q^*(z)=0$ so $z$ is not contained in $\Omega (\hat Q^*)$. This shows that every point $z\in (\hat F^*)^{-1}(\{\alpha\})$ is regular.  It now follows from Theorem \ref{theorem-minimal} and Theorem \ref{theorem-rational} that the fibre $(\hat F^*)^{-1}(\{\alpha\})$ over the non-zero $\alpha\in\cn^*$ is a complete {\it minimal submanifold} of $\cn^6_1$ of codimension two.
\end{example}

\section{The Complex Hyperbolic Space $\cn H^{2n-1}$}

In this section we construct complete non-holomorphic minimal submanifolds, of the odd-dimensional complex hyperbolic space $\cn H^{2n-1}$, of codimension two. Here we use the construction presented in Example \ref{example-hyperbolic-main}.
\medskip 

Define the $(4n-1)$-dimensional submanifold $\Sigma^{4n-1}$ of $\Omega^{2n}_1$ by 
$$
\Sigma^{4n-1}=\{z\in\Omega^{2n}_1\,|\, \ipp zz=-1\}.
$$ 
The map $\rho^*:S^1\times \Sigma^{4n-1}\to\Sigma^{4n-1}$ with $\rho^*:(e^{i\theta},z)\mapsto e^{i\theta} z$ is a smooth action of the Lie group $S^1$ on  $\Sigma^{4n-1}$ in $\Omega^{2n}_1$.  The quotient space of this action is the well-known complex hyperbolic space $\cn H^{2n-1}$.  The natural projection $\pi^*:\Sigma^{4n-1}\to\cn H^{2n-1}$ is a semi-Riemannian submersion with totally geodesic fibres and hence a harmonic morphism.

\begin{example}\label{example-complex-hyperbolic}
For $j,k\in\zn$ satisfying $1\le j\le n$ and  $n+1\le k\le 2n$, we define the function $\tilde\phi_{jk}^*:\Sigma^{4n-1}\to\cn$ with 
$$
\tilde\phi_{jk}^*:(z_1,\dots,z_{2n})\mapsto -\frac
{z_j\cdot\bar z_k}{\ipp zz}.
$$  
Then $\tilde\phi_{jk}^*:\Sigma^{4n-1}\to\cn$ is the restriction of the function
$\hat\phi_{jk}^*:\Omega^{2n}_1\to\cn$ to the set $\Sigma^{4n-1}$. An easy calculation shows that the tension field $\tau$ and the conformality operator $\kappa$ on $\Sigma^{4n-1}$ satisfy
$$
\tau(\tilde\phi_{jk}^*)=8n\cdot\tilde\phi_{jk}^*\ \ \text{and}\ \ \kappa(\tilde\phi_{jk}^*,\tilde\phi_{lm}^*)=4\cdot \tilde\phi_{jk}^*\,\tilde\phi_{lm}^*.
$$
This implies that the restrictions $\tilde P^*,\tilde Q^*:\Sigma^{4n-1}\to\cn$ of the polynomials $\hat P^*,\hat Q^*:\Omega^{2n}_1\to\cn$ to $\Sigma^{4n-1}$ form an eigenfamily.  Then we define the complex-valued harmonic morphism $\tilde F^*:\Omega (\tilde Q^*)\to\cn$ by $\tilde F^*(z)=\tilde P^*(z)/\tilde Q^*(z)$, where 
$$
\Omega(\tilde Q^*)=\{z\in \Sigma^{4n-1}\, |\, \tilde Q^*(z)\neq 0\}.
$$

The polynomial maps $\hat P^*,\hat Q^*:\Omega^{2n}_1\to\cn$ on $\Omega^{2n}_1$ induce  linearly independent functions $P^*,Q^*:\cn H^{2n-1}\to\cn$ on the complex hyperbolic space $\cn H^{2n-1}$, satisfying
$$
P^*:[z]\mapsto\sum_{j,k}a_{jk}\cdot\phi_{jk}^d\ \ \text{and}\ \ 
Q^*:[z]\mapsto\sum_{j,k}b_{jk}\cdot\phi_{jk}^d.
$$
Here the coefficients $a_{jk},b_{jk}\in\cn$ are chosen as explained in Example \ref{example-hyperbolic-main} and the induced functions $\phi_{jk}^*:\cn H^{2n-1}\to\cn$ satisfy 
$$
\phi_{jk}^*:[z]\mapsto -\frac
{z_j\cdot\bar z_k}{\ip zz}.
$$ 
The natural projection $\pi^*:\Sigma^{4n-1}\to\cn H^{2n-1}$ is a harmonic morphism so the complex valued functions $P^*,Q^*$ form an eigenfamily $\E_n^*$ on the complex hyperbolic space $\cn H^{2n-1}$.
\smallskip

Let $\Omega (Q^*)$ be the subset of $\cn H^{2n-1}$ given by 
$$
\Omega (Q^*)=\{[z]\in\cn H^{2n-1}\,|\,Q^*([z])\neq 0\}
$$
and define $F^*:\Omega(Q^*)\to\cn$ by $F^*([z])=P^*([z])/Q^*([z])$.  Then for a non-zero element  $\alpha\in\cn^*$ with $R(\alpha)\neq 0$, the fibre $(F^*)^{-1}(\{\alpha\})$ satisfies $(F^*)^{-1}(\{\alpha\})=\pi^* ((\tilde F^*)^{-1}(\{\alpha\}))$.  Since $\pi^*:\Sigma^{4n-1}\to\cn H^{2n-1}$ is a semi-Riemannian submersion the fibre is a {\it minimal submanifold} of codimension two satisfying 
$$
(F^*)^{-1}(\{\alpha\})=\{[z]\in\cn H^{2n-1}\,|\, P^*([z])=\alpha\cdot Q^*([z])\}.
$$
This is a {\it complete} and {\it non-holomorpic} submanifold of $\cn H^{2n-1}$.
\end{example}

\section{The Euclidean Space $\hn^{2n}\cong\cn^{2\times 2n}$}

Let $\cn^{2\times 2n}$ be the complex $4n$-dimensional vector space of complex 
$2\times 2n$ matrices.  We equip this space with its standard Euclidean metric $\ip \cdot\cdot :\cn^{2\times 2n}\times\cn^{2\times 2n}\to\rn$ satisfying 
$$\ip zw=\sum_{j,k}\Re (z_{jk}\cdot\bar w_{jk}).$$
Here the standard global coordinates on $\cn^{2\times 2n}$ are given by
$$z=
\begin{bmatrix}
	z_{11} & z_{12} & \cdots & z_{1,2n}  \\
	z_{21} & z_{22} & \cdots & z_{2,2n}
\end{bmatrix}.
$$

For two complex-valued functions $\hat\phi,\hat\psi:\cn^{2\times 2n}\to\cn$ the tension field $\tau$ and the conformal operator $\kappa$ are given by 
$$
\tau(\hat\phi)=4\cdot\sum_{j,k}\frac{\partial^2 \hat\phi}{\partial z_{jk}\partial\bar z_{jk}}
\ \ \text{and}\ \ \kappa(\hat\phi,\hat\psi)=2\cdot\sum_{j,k}\big(\frac{\partial\hat\phi}{\partial z_{jk}}\frac{\partial\hat\psi}{\partial\bar z_{jk}}+\frac{\partial\hat\phi}{\partial \bar z_{jk}}\frac{\partial\hat\psi}{\partial z_{jk}}\big).
$$

\begin{example}
\label{example-basic-complex-space}
For $j=1,2, \dots,n$ and $k=n+1,\dots ,2n$, we define the functions $\hat\phi_{jk}:\cn^{2\times 2n}\to\cn$ with 
$$
\hat\phi_{jk}:z\mapsto 
(z_{1j}\cdot\bar z_{1k}+z_{2j}\cdot\bar z_{2k}).$$ 
Then it is easily seen that we have a complex $n^2$-dimensional $(0,0)$-eigenfamily $\hat\E_m$ on $\cn^{2\times 2n}$ satisfying 
$$
\hat\E_n=\{\hat\phi_{jk}\,|\,j=1,2,\dots, n\ \ \text{and}\ \ k=n+1,\dots ,2n\}.
$$
\end{example}

\begin{example}\label{example-main}
Let $a,b\in\cn^9$ be two linearly independent elements satisfying
$$
a=(a_{14},a_{15},a_{16},a_{24},a_{25},a_{26},a_{34},a_{35},a_{36}),
$$
$$
b=(b_{14},b_{15},b_{16},b_{24},b_{25},b_{26},b_{34},b_{35},b_{36}).
$$
Further let the matrices $A,B\in\cn^{3\times 3}$ be given by 
$$
A=
\begin{bmatrix}
a_{14} & a_{15} & a_{16}  \\
a_{24} & a_{25} & a_{26}  \\
a_{34} & a_{35} & a_{36}
\end{bmatrix},\ \ 
B=
\begin{bmatrix}
b_{14} & b_{15} & b_{16}  \\
b_{24} & b_{25} & b_{26}  \\
b_{34} & b_{35} & b_{36}
\end{bmatrix}
$$
and assume that $\det B\neq 0$.  Then the third order polynomial $R:\cn\to\cn$ with $R(s)=\det (s\cdot B-A)$ is of the form 
$$
R_(s)=\det B\cdot s^3-c_2\cdot s^2+c_1\cdot s-\det A,
$$ 
for some complex numbers $c_1,c_2\in\cn$.
\smallskip
	
For $j=1,2,3$ and $k=4,5,6$ let us now define the two linearly independent polynomials $\hat P,\hat Q:\cn^6\to\cn$ by
$$
\hat P(z)=\sum_{j,k}a_{jk}\cdot (z_{1j}\cdot\bar z_{1k}+z_{2j}\cdot\bar z_{2k}),
$$
$$
\hat Q(z)=\sum_{j,k}b_{jk}\cdot (z_{1j}\cdot\bar z_{1k}+z_{2j}\cdot\bar z_{2k}).
$$
Further we define the complex-valued function $\hat F:\Omega (\hat Q)\to\cn$ by $\hat F(z)=\hat P(z)/\hat Q(z)$, where $\Omega(\hat Q)=\{z\in\cn^{12}\, |\, \hat Q(z)\neq 0\}$.
\smallskip
	
Let $\alpha\in\cn^*$ be a non-zero complex number such that $R(\alpha)\neq 0$.  Then a point $z\in \hat F^{-1}(\{\alpha\})$ in the fibre over $\alpha$ is critical if and only if 
$$
\frac{\partial\hat F}{\partial z_{11}}(z)
=\frac{\partial\hat F}{\partial z_{12}}(z)
=\frac{\partial\hat F}{\partial z_{13}}(z)=0
=\frac{\partial\hat F}{\partial z_{21}}(z)
=\frac{\partial\hat F}{\partial z_{22}}(z)
=\frac{\partial\hat F}{\partial z_{23}}(z)
$$
and
$$
\frac{\partial\hat F}{\partial \bar z_{14}}(z)
=\frac{\partial\hat F}{\partial \bar z_{15}}(z)
=\frac{\partial\hat F}{\partial \bar z_{16}}(z)
=0
=\frac{\partial\hat F}{\partial \bar z_{24}}(z)
=\frac{\partial\hat F}{\partial \bar z_{25}}(z)
=\frac{\partial\hat F}{\partial \bar z_{26}}(z).
$$
It is easily seen that this is equivalent to the following systems of equations
$$
\Big(t\cdot B - A\Big)\cdot 
\begin{bmatrix}
\bar z_{41} & \bar z_{42} & z_{11} & z_{12}\\
\bar z_{51} & \bar z_{52} & z_{21} & z_{22}\\
\bar z_{61} & \bar z_{62} & z_{31} & z_{32}
\end{bmatrix}=0.
$$

Since $R(\alpha)\neq 0$ this implies that $z=0$ which is not contained in $\Omega (\hat Q)$. This shows that every point $z\in \hat F^{-1}(\{\alpha\})$ is regular.  It now follows from Theorem \ref{theorem-minimal} and Theorem \ref{theorem-rational} that the fibre $\hat F^{-1}(\{\alpha\})$ over the non-zero $\alpha\in\cn^*$ is a complete {\it minimal submanifold} of $\cn^{12}$ of codimension two.
\end{example}
\smallskip

\section{The Quaternionic Projective Space $\hn P^{2n-1}$}

In this section we provide complete minimal submanifolds, of the quaternionic projective space $\hn P^{2n-1}$, of codimension two.
\medskip 

The map $\rho:S^3\times S^{8n-1}\to S^{8n-1}$ with $\rho:(q,z)\mapsto z\cdot q$ is a smooth action of the Lie group $S^3$, of unit quaternions, on the unit sphere $S^{8n-1}$ in $\cn^{2\times 2n}$. The quotient space of this action is the well-known  quaternionic projective space $\hn P^{2n-1}$.  The natural projection $\pi:S^{8n-1}\to\hn P^{2n-1}$ is a Riemannian submersion with totally geodesic fibres and hence a harmonic morphism.

\begin{example}\label{example-quaternionic-projective}
For $j=1,2,\dots,n$ and $k=n+1,\dots, 2n$, we define the $S^3$-invariant functions $\tilde\phi_{jk}:S^{8n-1}\to\cn$ on the unit sphere in $S^{8n-1}$ in $\cn^{2\times 2n}$ by 
$$
\tilde\phi_{jk}:z\mapsto \frac
{(z_{1j}\cdot\bar z_{1k}+z_{2j}\cdot\bar z_{2k})}{\ip zz}.
$$  
Then $\tilde\phi_{jk}:S^{8n-1}\to\cn$ is the restriction of 
the function $\hat\phi_{jk}:\cn^{2\times 2n}\to\cn$ to the unit sphere $S^{8n-1}$. An easy calculation shows that the tension field $\tau$ and the conformality operator $\kappa$ on $S^{8n-1}$ satisfy
$$
\tau(\tilde\phi_{jk})=-\,8n\cdot\tilde\phi_{jk}\ \ \text{and}\ \ \kappa(\tilde\phi_{jk},\tilde\phi_{lm})=-\,4\cdot \tilde\phi_{jk}\,\tilde\phi_{lm}.
$$
This shows that we have a complex $n^2$-dimensional eigenfamily $\tilde\E_n$ on the unit sphere $S^{8n-1}$.	 This implies that the restrictions $\tilde P,\tilde Q:S^{8n-1}\to\cn$ of the polynomials $\hat P,\hat Q:\cn^{2\times 2n}\to\cn$ to the unit sphere $S^{8n-1}$ form an eigenfamily.  Then we define the complex-valued harmonic morphism $\tilde F:\Omega (\tilde Q)\to\cn$ by $\tilde F(z)=\tilde P(z)/\tilde Q(z)$, where 
$$
\Omega(\tilde Q)=\{z\in S^{8n-1}\, |\, \tilde Q(z)\neq 0\}.
$$

The map $\tilde F:\Omega (\tilde Q)\to\cn$ is invariant under the $S^3$-action on the unit sphere $S^{8n-1}$ i.e. $\hat F(z\cdot q)=\hat F(z)$ for all $q\in S^3$.

For $j=1,2,\dots, n$ and  $k=n+1,\dots ,2n$, we define the functions  $\phi_{jk}:\hn P^{2n-1}\to\cn$ by	
$$
\phi_{jk}:[z]\mapsto \frac
{(z_{1j}\cdot\bar z_{1k}+z_{2j}\cdot\bar z_{2k})}{\ip zz}.
$$  
The natural projection $\pi:S^{8n-1}\to\hn P^{n-1}$ is a harmonic Riemannian submersion so it follows from Corollary 3.5 of \cite{Geg-Gud-1} that the tension field $\tau$ and the conformality operator $\kappa$ on $\hn P^{2n-1}$ satisfy
$$\tau(\phi_{jk})=-\,8n\cdot\phi_{jk}\ \ \text{and}\ \ \kappa(\phi_{jk},\phi_{lm})=-\,4\cdot \phi_{jk}\,\phi_{lm}.$$
The polynomial maps $\hat P,\hat Q:\cn^{2\times 2n}\to\cn$ on $\hn^{2n}$ induce  linearly independent functions $P,Q:\hn P^{2n-1}\to\cn$ on the complex projective space $\hn P^{2n-1}$, satisfying
$$
P:[z]\mapsto\sum_{j,k}a_{jk}\cdot\phi_{jk}\ \ \text{and}\ \ 
Q:[z]\mapsto\sum_{j,k}b_{jk}\cdot\phi_{jk}.
$$
The induced maps $P,Q$ form an eigenfamily $\E_n$ on the quaternionic projective space $\hn P^{2n-1}$.
	 
Let $\Omega (Q)$ be the subset of $\hn P^{2n-1}$ given by 
$$
\Omega (Q)=\{[z]\in\hn P^{2n-1}\,|\,Q([z])\neq 0\}
$$
and define $F:\Omega(Q)\to\cn$ by $F([z])=P([z])/Q([z])$.  Then for a non-zero element  $\alpha\in\cn^*$ with $R(\alpha)\neq 0$, the fibre $F^{-1}(\{\alpha\})$ satisfies $F^{-1}(\{\alpha\})=\pi (\tilde F^{-1}(\{\alpha\}))$.  Since $\pi:S^{8n-1}\to\hn P^{2n-1}$ is a Riemannian submersion the fibre is a {\it minimal submanifold} of codimension two, satisfying 
$$
F^{-1}(\{\alpha\})=\{[z]\in\hn P^{2n-1}\,|\, P([z])=\alpha\cdot Q([z])\}.
$$
This is a {\it complete} submanifold of $\hn P^{2n-1}$.
\end{example}

\section{The Semi-Euclidean Space $\hn^{2n}_1=\cn^{2\times 2n}_1$}

Let $\cn^{2\times 2n}_1$ be the complex $4n$-dimensional vector space equipped with its standard semi-Euclidean metric $\ipp \cdot\cdot :\cn^{2\times 2n}_1\times\cn^{2\times 2n}_1\to\rn$ satisfying 
$$
\ipp zw=\Re \Big(-(z_{11}\cdot \bar w_{11}+z_{21}\cdot \bar w_{21}) +\sum_{k=2}^{2n}(z_{1k}\cdot\bar w_{1k}+z_{2k}\cdot\bar w_{2k})\Big).
$$
Let $\Lambda^{4n}_1$ be the open subset $\Lambda^{4n}_1=\{z\in\cn^{2\times 2n}_1\,|\, \ipp zz<0\}$ of $\cn^{2\times 2n}_1$.

For two complex-valued functions $\hat\phi^*,\hat\psi^*:\Lambda^{4n}_1\to\cn$ the tension field $\tau$ and the conformal operator $\kappa$ on $\Lambda^{4n}_1$, are given by 
$$
\tau(\hat\phi^*)=4\cdot\Big(-\big(\frac{\partial^2 \hat\phi^*}{\partial z_{11}\partial\bar z_{11}}
+\frac{\partial^2 \hat\phi^*}{\partial z_{21}\partial\bar z_{21}}\big)
+ \sum_{k=2}^{2n}\big(\frac{\partial^2 \hat\phi^*}{\partial z_{1k}\partial\bar z_{1k}}+\frac{\partial^2 \hat\phi^*}{\partial z_{2k}\partial\bar z_{2k}}\big) \Big)
$$
and
\begin{eqnarray*}
& &\kappa(\hat\phi^*,\hat\psi^*)\\
&=&2\cdot\Big(-\big(\frac{\partial\hat\phi^*}{\partial z_{11}}\frac{\partial\hat\psi^*}{\partial\bar z_{11}}+\frac{\partial\hat\phi^*}{\partial \bar z_{11}}\frac{\partial\hat\psi^*}{\partial z_{11}}
+\frac{\partial\hat\phi^*}{\partial z_{21}}\frac{\partial\hat\psi^*}{\partial\bar z_{21}}+\frac{\partial\hat\phi^*}{\partial \bar z_{21}}\frac{\partial\hat\psi^*}{\partial z_{21}}\big)
\\
& &+\sum_{k=2}^{2n}\big(\frac{\partial\hat\phi^*}{\partial z_{1k}}\frac{\partial\hat\psi^*}{\partial\bar z_{1k}}+\frac{\partial\hat\phi^*}{\partial \bar z_{1k}}\frac{\partial\hat\psi^*}{\partial z_{1k}}
+\frac{\partial\hat\phi^*}{\partial z_{2k}}\frac{\partial\hat\psi^*}{\partial\bar z_{2k}}+\frac{\partial\hat\phi^*}{\partial \bar z_{2k}}\frac{\partial\hat\psi^*}{\partial z_{2k}}
\big)\Big).
\end{eqnarray*}
Here the standard global coordinates on $\Lambda^{4n}_1$ are given by
$$z=
\begin{bmatrix}
z_{11} & z_{12} & \cdots & z_{1,2n}  \\
z_{21} & z_{22} & \cdots & z_{2,2n}
\end{bmatrix}.
$$

\begin{example}
\label{example-basic-complex-space}
For $j,k\in\zn$ satisfying $1\le j\le n$ and  $n+1\le k\le 2n$, define the functions  $\hat\phi^*_{jk}:\Lambda^{4n}_1\to\cn$ with 
$$
\hat\phi^*_{jk}:z\mapsto ({z_{1j}\cdot\bar z_{1k}+z_{2j}\cdot\bar z_{2k}}).
$$ 
Then we have a complex $n^2$-dimensional $(0,0)$-eigenfamily $\hat\E_n^*$ on $\Lambda^{4n}_1$ satisfying 
$$
\hat\E_n^*=\{\hat\phi_{jk}^*\,|\,1\le j\le n\ \ \text{and}\ \ n+1\le k\le 2n\}.
$$
\end{example}

\begin{example} 
Let $a,b\in\cn^9$ be two linearly independent elements, satisfying
$$
a=(a_{14},a_{15},a_{16},a_{24},a_{25},a_{26},a_{34},a_{35},a_{36}),
$$
$$ b=(b_{14},b_{15},b_{16},b_{24},b_{25},b_{26},b_{34},b_{35},b_{36}).
$$
	Further let the matrices $A,B\in\cn^{3\times 3}$ be given by 
	$$A=
	\begin{bmatrix}
		a_{14} & a_{15} & a_{16}  \\
		a_{24} & a_{25} & a_{26}  \\
		a_{34} & a_{35} & a_{36}
	\end{bmatrix},\ \ B=
	\begin{bmatrix}
		b_{14} & b_{15} & b_{16}  \\
		b_{24} & b_{25} & b_{26}  \\
		b_{34} & b_{35} & b_{36}
	\end{bmatrix}
	$$
	and assume that $\det B\neq 0$.  Then the third order polynomial $R:\cn\to\cn$ with $R(s)=\det (s\cdot B-A)$ is of the form 
	$$
	R(s)=\det B\cdot s^3-c_2\cdot s^2+c_1\cdot s-\det A,
	$$ 
	for some complex numbers $c_1,c_2\in\cn$.
	\smallskip
	
	For $j=1,2,3$ and $k=4,5,6$ let us now define the two linearly independent polynomials $\hat P^*,\hat Q^*:\Lambda^{12}_1\to\cn$ by
	$$
	\hat P^*(z)=\sum_{j,k}a_{jk}\cdot
	\hat\phi_{jk}^*
	\ \ \text{and}\ \ 
	\hat Q^*(z)=\sum_{j,k}b_{jk}\cdot
	\hat\phi_{jk}^*
	$$
	Further we define the complex-valued function $\hat F^*:\Omega (\hat Q^*)\to\cn$ by $\hat F^*(z)=\hat P^*(z)/\hat Q^*(z)$, where $\Omega(\hat Q^*)=\{z\in\Lambda^{12}_1\, |\, \hat Q^*(z)\neq 0\}$.
	\smallskip
	
	Let $\alpha\in\cn^*$ be a non-zero complex number such that $R(\alpha)\neq 0$.  Then a point $z\in (\hat F^*)^{-1}(\{\alpha\})$ in the fibre over $\alpha$ is critical if and only if 
	$$
	 \frac{\partial\hat F^*}{\partial z_{11}}(z)
	=\frac{\partial\hat F^*}{\partial z_{12}}(z)
	=\frac{\partial\hat F^*}{\partial z_{13}}(z)=0
	=\frac{\partial\hat F^*}{\partial z_{21}}(z)
	=\frac{\partial\hat F^*}{\partial z_{22}}(z)
	=\frac{\partial\hat F^*}{\partial z_{23}}(z)
	$$
	and
	$$
	 \frac{\partial\hat F^*}{\partial \bar z_{14}}(z)
	=\frac{\partial\hat F^*}{\partial \bar z_{15}}(z)
	=\frac{\partial\hat F^*}{\partial \bar z_{16}}(z)=0
	=\frac{\partial\hat F^*}{\partial \bar z_{24}}(z)
	=\frac{\partial\hat F^*}{\partial \bar z_{25}}(z)
	=\frac{\partial\hat F^*}{\partial \bar z_{26}}(z).
	$$
	It is easily seen that this is equivalent to the following systems of equations
	$$
	\Big(t\cdot B - A\Big)\cdot 
	\begin{bmatrix}
		\bar z_{41} & \bar z_{42} & z_{11} & z_{12}\\
		\bar z_{51} & \bar z_{52} & z_{21} & z_{22}\\
		\bar z_{61} & \bar z_{62} & z_{31} & z_{32}
	\end{bmatrix}=0.
	$$
	Since $R(\alpha)\neq 0$ this implies that $\hat Q^*(z)=0$ so $z$ is not contained in $\Omega (\hat Q^*)$. This shows that every point $z\in (\hat F^*)^{-1}(\{\alpha\})$ is regular.  It now follows from Theorem \ref{theorem-minimal} and Theorem \ref{theorem-rational} that the fibre $(\hat F^*)^{-1}(\{\alpha\})$ over the non-zero $\alpha\in\cn^*$ is a complete {\it minimal submanifold} of $\cn^{12}_1$ of codimension two.
\end{example}

\section{The Quaterionic Hyperbolic Space $\hn H^{2n-1}$}

In this section we provide complete minimal submanifolds, of the quaternionic hyperbolic space $\hn H^{2n-1}$, of codimension two.
\medskip 

Define the $(8n-1)$-dimensional submanifold $\Sigma^{8n-1}$ of $\Lambda^{4n}_1$ by 
$$
\Sigma^{8n-1}=\{z\in\Lambda^{4n}_1\,|\, \ipp zz=-1\}.
$$ 
The map $\rho:S^3\times \Sigma^{8n-1}\to \Sigma^{8n-1}$ with $\rho:(q,z)\mapsto z\cdot q$ is a smooth action of the Lie group $S^3$ on $\Sigma^{8n-1}$ in $\Lambda^{4n}_1$. The quotient space of this action is the quaternionic hyperbolic space $\hn H^{2n-1}$.  The natural projection $\pi^*:\Sigma^{8n-1}\to\hn H^{2n-1}$ is a well-known semi-Riemannian submersion with totally geodesic fibres and hence a harmonic morphism.

\begin{example}\label{example-quaternionic-hyperbolic}
For $j,k\in\zn$ satisfying $1\le j\le n$ and  $n+1\le k\le 2n$, we define the function $\tilde\phi_{jk}^*:\Sigma^{8n-1}\to\cn$ with 
$$
\tilde\phi_{jk}^*:z\mapsto -\frac
{(z_{1j}\cdot\bar z_{1k}+z_{2j}\cdot\bar z_{2k})}{\ipp zz}.
$$  
Then $\tilde\phi_{jk}^*:\Sigma^{8n-1}\to\cn$ is clearly the restriction of the function $\hat\phi_{jk}^*:\Lambda^{4n}_1\to\cn$ to the set $\Sigma^{8n-1}$. An easy calculation shows that the tension field $\tau$ and the conformality operator $\kappa$ on $\Sigma^{8n-1}$ satisfy
$$
\tau(\tilde\phi_{jk}^*)=8n\cdot\tilde\phi_{jk}^*\ \ \text{and}\ \ \kappa(\tilde\phi_{jk}^*,\tilde\phi_{lm}^*)=4\cdot \tilde\phi_{jk}^*\,\tilde\phi_{lm}^*.
$$
This shows that we have a complex $n^2$-dimensional eigenfamily $\tilde\E_n^*$ on the set $\Sigma^{8n-1}$.	 This implies that the restrictions $\tilde P^*,\tilde Q^*:\Sigma^{8n-1}\to\cn$ of the polynomials $\hat P^*,\hat Q^*:\Lambda^{4n}_1\to\cn$ to the set $\Sigma^{8n-1}$ form an eigenfamily.  Then we define the complex-valued function $\tilde F^*:\Omega (\tilde Q^*)\to\cn$ by $\tilde F^*(z)=\tilde P^*(z)/\tilde Q^*(z)$, where 
$$
\Omega(\tilde Q^*)=\{z\in \Sigma^{8n-1}\, |\, \tilde Q^*(z)\neq 0\}.
$$
For $j,k\in\zn$ satisfying $1\le j\le n$ and  $n+1\le k\le 2n$, we define the functions $\phi_{jk}^*:\hn H^{2n-1}\to\cn$ by	
$$
\phi_{jk}^*:[z_{11},\dots,z_{2,2n}]\mapsto -\frac
{(z_{1j}\cdot\bar z_{1k}+z_{2j}\cdot\bar z_{2k})}{\ipp zz}.
$$  
The projection map $\pi^*:\Sigma^{8n-1}\to\hn H^{2n-1}$ is a semi-Riemannian submersion so the tension field $\tau$ and the conformality operator $\kappa$ on $\hn H^{2n-1}$ satisfy
$$
\tau(\phi_{jk}^*)=8n\cdot\phi_{jk}^*\ \ \text{and}\ \ \kappa(\phi_{jk}^*,\phi_{lm}^*)=4\cdot \phi_{jk}^*\,\phi_{lm}^*.
$$
The polynomial maps $\hat P^*,\hat  Q^*:\Lambda^{4n}_1\to\cn$ on $\Lambda^{4n}_1$ induce  linearly independent functions $P^*,Q^*:\hn H^{2n-1}\to\cn$ on the complex hyperbolic space $\hn H^{2n-1}$, satisfying
$$
P^*:[z]\mapsto\sum_{j,k}a_{jk}\cdot\phi_{jk}^*\ \ \text{and}\ \ 
Q^*:[z]\mapsto\sum_{j,k}b_{jk}\cdot\phi_{jk}^*.
$$
Here the coefficients $a_{jk},b_{jk}\in\cn$ are chosen precisely the same way as explained in Example \ref{example-quaternionic-hyperbolic}.  The induced maps $P^*,Q^*$ form an eigenfamily $\E_n^*$ on the quaternionic  hyperbolic space $\hn H^{2n-1}$.
	\smallskip 

Let $\Omega (Q^*)$ be the subset of $\hn H^{2n-1}$ given by 
$$
\Omega (Q^*)=\{[z]\in\hn H^{2n-1}\,|\,Q^*([z])\neq 0\}
$$
and define $F^*:\Omega(Q^*)\to\cn$ by $F^*([z])=P^*([z])/Q^*([z])$.  Then for a non-zero element  $\alpha\in\cn^*$ with $R(\alpha)\neq 0$, the fibre $(F^*)^{-1}(\{\alpha\})$ satisfies $(F^*)^{-1}(\{\alpha\})=\pi^* ((\tilde F^*)^{-1}(\{\alpha\}))$.  Since $\pi^*:\Sigma^{8n-1}\to\hn H^{2n-1}$ is a semi-Riemannian submersion the fibre is a {\it minimal submanifold} of codimension two satisfying 
$$
(F^*)^{-1}(\{\alpha\})=\{[z]\in\hn H^{2n-1}\,|\, P^*([z])=\alpha\cdot Q^*([z])\}.
$$
This is clearly a {\it complete} submanifold of $\hn H^{2n-1}$.
\end{example}

\section{Acknowledgements}

The author is grateful to Thomas Jack Munn for useful discussions on this work.



\begin{thebibliography}{99}
	
	
\bibitem{Bai-Eel}
P. Baird, J. Eells,
{\it A conservation law for harmonic maps},
Geometry Symposium Utrecht 1980, Lecture Notes in Mathematics {\bf 894}, Springer (1981), 1-25.
	
	
\bibitem{Bai-Woo-book} 
P. Baird and J. C. Wood,
{\it Harmonic Morphisms Between Riemannian Manifolds},
London Math. Soc. Monogr. {\bf 29},	
Oxford Univ. Press (2003).
	
	
\bibitem{Bur-Raw}
F. E. Burstall, J. H. Rawnsley,
{\it Twistor Theory for Riemannian Symmetric Spaces.
With Applications to Harmonic Maps of Riemann Surfaces},
Lecture Notes in Math. {\bf 1424}, Springer (1990).
	
	
\bibitem{Bur-Gue}
F.E. Burstall, M.A. Guest, 
{\it Harmonic two-spheres in compact symmetric spaces, revisited},
Math. Ann. {\bf 309} (1997), 541-752.
	
	
\bibitem{Cal}
E. Calabi, 
{\it Minimal immersions of surfaces in Euclidean spheres}, 
J. Diff. Geom. {\bf 1} (1967), 111-125.
	
	
\bibitem{Eel-Sam}
J. Eells, J. H. Sampson,
{\it Harmonic mappings of Riemannian manifolds},
Amer. J. Math. {\bf 86} (1964), 109-160.
	
	
\bibitem{Eel-Woo}
J. Eells, J.C. Wood, 
{\it Harmonic maps from surfaces to complex projective spaces}, 
Adv. in Math. {\bf 49}, (1983), 217-263.


\bibitem{Fug-2} 
B. Fuglede,
{\it Harmonic morphisms between semi-riemannian manifolds}, 
Ann. Acad. Sci. Fennicae {\bf 21} (1996), 31-50.


\bibitem{Geg-Gud-1}
J. M. Gegenfurtner, S. Gudmundsson,
{\it Compact minimal submanifolds of the Riemannian symmetric spaces $\SU n/\SO n$, $\Sp n/\U n$, $\SO{2n}/\U n$, $\SU{2n}/\Sp n$ via complex-valued eigenfunctions},  
Ann. Global Anal. Geom. {\bf 66} (2024), 13.


\bibitem{Gha-Gud-5}
E. Ghandour, S. Gudmundsson,
{\it Explicit harmonic morphisms and $p$-harmonic functions from the complex and quaternionic Grassmannians}, 
Ann. Global Anal. Geom. {\bf 64} (2023), 15.

	
\bibitem{Gud-bib}
S. Gudmundsson,
{\it The Bibliography of Harmonic Morphisms},
{\tt www.matematik.lu.se/ matematiklu/personal/sigma/harmonic/bibliography.html}


\bibitem{Gud-8}
S. Gudmundsson,
{\it On the existence of harmonic morphisms from symmetric spaces of rank one},
Manuscripta Math. {\bf 93} (1997), 421-433.


\bibitem{Gud-Mon-Rat-1}
S. Gudmundsson, S. Montaldo, A. Ratto,
{\it Biharmonic functions on the classical compact simple Lie groups},
J. Geom. Anal. {\bf 28} (2018), 1525-1547.
	
	
\bibitem{Gud-Mun-1}
S. Gudmundsson, T. J. Munn,
{\it Minimal submanifolds via complex-valued eigenfunctions},
J. Geom. Anal. {\bf 34} (2024), 190.
	
	
\bibitem{Gud-Sak-1}	
S. Gudmundsson, A. Sakovich,
{\it Harmonic morphisms from the classical compact semisimple Lie groups},
Ann. Global Anal. Geom. {\bf 33} (2008), 343-356.

	
\bibitem{Gud-Sve-1}
S. Gudmundsson and M. Svensson,
{\it Harmonic morphisms from the Grassmannians and their non-compact duals},
Ann. Global Anal. Geom. {\bf 30} (2006), 313-333.	
	
\bibitem{Hel}
S. Helgason, 
{\it Differential Geometry, Lie Groups, and Symmetric Spaces}, Grad. Stud. Math. {\bf 34}, American Mathematical Society, (2001).

	
\bibitem{Uhl}
K. Uhlenbeck, 
{\it Harmonic maps into Lie groups: classical solutions of the chiral model},
J. Differential Geom. {\bf 30} (1989), 1–50.
	
\end{thebibliography}
\end{document}